\documentclass{amsart}

\usepackage{amssymb,amsfonts,latexsym}

\numberwithin{equation}{section}

\newtheorem{theorem}{Theorem}[section]

\newtheorem{proposition}[theorem]{Proposition}

\theoremstyle{definition}

\newtheorem *{Theorem A}{Theorem A}
\newtheorem *{Theorem B}{Theorem B}
\newtheorem *{Corollary C}{Corollary C}

\newcommand{\btu}{\bar{\tau}}

\newcommand{\tut}{\eta_t}

\newcommand{\Cn}{C_{2^{n-1}}}

\newcommand{\hbt}{{\widehat {\beta}}}

\newcommand{\la}{\langle\,}

\newcommand{\ra}{\,\rangle}

\newcommand{\ut}{\text{Aut}}

\newcommand{\ben}{\begin{enumerate}}

\newcommand{\een}{\end{enumerate}}

\newcommand{\bx}{{\bar x}}

\newcommand{\C}{{\mathbb C}}

\begin{document}

\title{ A separable deformation of the quaternion group algebra}

\author{Nurit Barnea}

\author{Yuval Ginosar}

\address{Department of Mathematics, University of Haifa, Haifa 31905, Israel}
\email{ginosar@math.haifa.ac.il}

\date{\today}

\begin{abstract}
The Donald-Flanigan conjecture asserts that for any finite group $G$ and any field $k$,
the group algebra $kG$ can be deformed to a separable algebra.
The minimal unsolved instance, namely the quaternion group $Q_8$
over a field $k$ of characteristic 2 was considered as a counterexample.
We present here a separable deformation of $kQ_8$.
In a sense, the conjecture for any finite group is open again.
\end{abstract}

\maketitle
\section{Introduction}
In their paper \cite{DF}, J.D. Donald and F.J. Flanigan conjectured that any group algebra
$kG$ of a finite group $G$ over a field $k$ can be deformed to a semisimple algebra
even in the modular case, namely where the order of $G$ is not invertible in $k$.
A more customary formulation of the Donald-Flanigan (DF) conjecture is by demanding
that the deformed algebra $[kG]_t$ should be separable, i.e. it remains semisimple when tensored
with the algebraic closure of its base field.
If, additionally, the dimensions of the simple components
of $[kG]_t$ are in one-to-one correspondence with those of the complex
group algebra $\C G$, then $[kG]_t$ is called a {\it strong} solution to the problem.

The DF conjecture was solved for groups $G$ which have either a cyclic $p$-Sylow subgroup over an
algebraically closed field
\cite{S} or a normal abelian $p$-Sylow subgroup \cite{GS96}
where $p=$char$(k)$, and for all but six reflection groups in any
characteristic \cite{GS97,GGS,PS}.
In \cite{GG}, it is claimed that the group algebra $kQ_8$, where
$$Q_8=\la\sigma,\tau| \sigma^4=1, \tau\sigma=\sigma^3\tau, \sigma^2=\tau^2\ra  $$
is the quaternion group of order 8 and $k$ a field of characteristic
$2$, does not admit a separable deformation. This result allegedly
gives a counterexample to the DF conjecture. However, as observed by M. Schaps,
the proof apparently contains an error (see \S\ref{ack}).

The aim of this note is to present a separable deformation of $kQ_8$, where $k$ is any field
of characteristic 2,
reopening the DF conjecture.

\section{Preliminaries}
Let $k[[t]]$ be the ring of formal power series over $k$, and let $k((t))$ be its field of fractions.
Recall that the deformed algebra $[kG]_t$ has the same underlying $k((t))$-vector space as
$k((t))\otimes_k kG$,
with multiplication defined on basis elements
\begin{equation}\label{prod}
g_1* g_2:=g_1g_2+\sum_{i\geq 1} \Psi_{i}(g_1,g_2)t^{i},\ \ g_1,g_2\in G
\end{equation}
and extended $k((t))$-linearly (such that $t$ is central).
Here $g_1g_2$ is the group multiplication.
The functions $\Psi_{i}:G\times G\to kG$ satisfy certain
cohomological conditions induced by the associativity of
$[kG]_t$ \cite[\S 1 ; \S 2]{G}.

Note that the set of equations (\ref{prod})
determines a multiplication on the free $k[[t]]$-module $\Lambda_t$
spanned by the elements $\{g\}_{g\in G}$
such that $kG\simeq \Lambda_t/\langle t\Lambda_t \rangle$
and $[kG]_t\simeq\Lambda_t\otimes_{k[[t]]}k((t))$.
In a more general context, namely over a domain $R$ which is not necessarily local,
the $R$-module $\Lambda_t$ which determines the deformation, is required only to be {\it flat}
rather than free \cite[\S 1]{ES}.

In what follows, we shall define the deformed algebra $[kG]_t$ by using
generators and relations. These will implicitly determine the set of equations (\ref{prod}).


\section{Sketch of the construction}

Consider the extension
\begin{equation}
[\beta]: 1\to C_4\to Q_8\to C_2\to 1,
\end{equation}
where $C_2=\la\bar\tau\ra  $ acts on $C_4=\la\sigma\ra  $ by
$$\begin{array}{rccl}
\eta:&C_2&\to &\ut (C_4)\\
\eta(\btu):&\sigma&\mapsto &\sigma^3(=\sigma^{-1}),
\end{array}$$
and the associated 2-cocycle
$\beta:C_2\times C_2\to C_4$ is given by
$$\beta(1,1)=\beta(1,\btu)=\beta(\btu,1)=1,\beta(\btu,\btu)=\sigma^2.$$
The group algebra $kQ_8$ ($k$ any field) is isomorphic to
the quotient $kC_4[y;\eta]/\la q(y)\ra  $,
where $kC_4[y;\eta]$ is a skew polynomial ring \cite[\S 1.2]{MR},
whose indeterminate $y$ acts on the ring of coefficients $kC_4$ via
the automorphism $\eta(\btu)$ (extended linearly) and where
\begin{equation}\label{q}
q(y):=y^2-\sigma^2  \in kC_4[y;\eta]
\end{equation}
is central. The above isomorphism is established by identifying $\tau$
with the indeterminate $y$.

Suppose now that Char$(k)=2$. The deformed algebra $[k{Q_8}]_t$ is constructed as follows.

In
\S \ref{C4.1} the subgroup algebra $kC_4$ is deformed to a
separable algebra $[k{C_4}]_t$ which is isomorphic to
$K\oplus k((t))\oplus k((t))$, where $K$ is a separable field extension
of $k((t))$ of degree 2.

The next step (\S\ref{C4.2}) is to
construct an automorphism $\tut$ of $[k{C_4}]_t$ which agrees
with the action of $C_2$ on $kC_4$ when specializing  $t=0$.
This action fixes all three primitive idempotents of $[{kC_4}]_t$.
By that we obtain the skew polynomial ring $[k{C_4}]_t[y;\tut]$.

In \S \ref{hbt} we deform $q(y)=y^2+\sigma^2$ to $q_t(y)$,
a separable polynomial of degree 2 in the center of $[kC_4]_t[y;\tut]$.

By factoring out the two-sided ideal generated by $q_t(y)$, we establish the deformation
$$[{kQ_8}]_t:=[{kC_4}]_t[y;\tut]/\la q_t(y)\ra  .$$

In \S\ref{sep} we show that $[{kQ_8}]_t$ as above is separable.
Moreover, passing to the algebraic closure $\overline{k((t))}$
we have
$$[kQ_8]_t\otimes_{k((t))} \overline{k((t))}\simeq \bigoplus_{i=1}^4\overline{k((t))}\oplus M_2(\overline{k((t))}).$$
This is a strong solution to the DF conjecture since
its decomposition to simple components is the same as
$$\C Q_8\simeq\bigoplus\limits_{i=1}^4\C\oplus M_2(\C).$$
\section{A Deformation of $kC_4[y;\eta]$}\label{C4}

\subsection{}  \label{C4.1}
We begin by constructing $[{kC_4}]_t$,  $C_4=\la \sigma\ra$.
Recall that $$kC_4\simeq k[x]/\la x^4+1\ra$$
by identifying $\sigma$ with $x+\la x^4+1\ra$. We deform the polynomial $x^4+1$ to a separable
polynomial $p_t(x)$ as follows.

Let $k[[t]]^*$ be the group of invertible elements of $k[[t]]$ and denote by
$$U:=\{1+zt| z\in k[[t]]^*\}$$ its subgroup of 1-units (when $k=\mathbb{F}_2$,
$U$ is equal to $k[[t]]^*$).

Let
$$a \in k[[t]]\setminus k[[t]]^*$$
be a non-zero element, and let
$$b,c,d\in U, (c\neq d),$$
such that
$$\pi(x):=x^2+ax+b$$
is an irreducible (separable) polynomial in $k((t))[x]$.
Let
$$p_t(x):=\pi(x)(x+c)(x+d)\in k((t))[x].$$
Then the quotient $k((t))[x]/\la p_t(x)\ra  $ is
isomorphic to the direct sum $K\oplus k((t))\oplus k((t))$, where
$K:=k((t))[x]/\la \pi(x)\ra  $. The field extension $K/k((t))$ is separable and of dimension 2.

Note that $p_{t=0}(x)=x^4+1$ and that only lower order terms of the polynomial were deformed.
Hence, the quotient $k[[t]][x]/\la p_t(x)\ra  $ is $k[[t]]$-free and
$k((t))[x]/\la p_t(x)\ra  $
indeed defines a deformation $[{kC_4}]_t$ of $kC_4\simeq k[x]/\la x^4+1\ra  $.
The new multiplication $\sigma^i*\sigma^j$ of basis elements
(\ref{prod}) is determined by identifying
$\sigma^i$ with $\bx^i:=x^i+\la p_t(x)\ra$.
We shall continue to use the term $\bx$ in $[{kC_4}]_t$ rather than $\sigma$.

Assume further that there exists $w\in k[[t]]$ such that
\begin{equation}\label{assume}
(x+w)(x+c)(x+d)=x\pi(x)+a
\end{equation}
(see example \ref{example}).
Then $K\simeq([{kC_4}]_t)e_1$, where
\begin{equation} \label{idem1}
e_1=\frac{(\bx+w)(\bx+c)(\bx+d)}{a}  .
\end{equation}
The two other primitive idempotents of $[{kC_4}]_t$ are
\begin{equation}\label{idem23}
e_2=\frac{c(\bx+d)\pi(\bx)}{a(c+d)}, \ \ e_3=\frac{d(\bx+c)\pi(\bx)}{a(c+d)}.
\end{equation}

\subsection{} \label{C4.2}

Let $$\tut:k((t))[x]\rightarrow k((t))[x]$$
be an algebra endomorphism determined by its value on the
generator $x$ as follows.
\begin{equation}\label{p8}
\eta_t(x):=x\pi(x)+x+a.
\end{equation}
We compute $\tut(\pi(x))$,  $\tut(x+c)$ and $\tut(x+d)$:
\begin{equation*}
\begin{split}
\tut(\pi(x))=\tut(x)^2+a\tut(x)+b&=x^2\pi(x)^2+x^2+a^2+ax\pi(x)+ax+a^2+b\\
&=\pi(x)(x^2\pi(x)+ax+1).
\end{split}
\end{equation*}
By (\ref{assume}),
\begin{equation}\label{p1}
\tut(\pi(x))=\pi(x)+x(x+w)p_t(x)\in\la \pi(x)\ra  .
\end{equation}
Next,
$$\tut(x+c)=x\pi(x)+x+a+c.$$
By (\ref{assume}),
\begin{equation}\label{p2}
\tut(x+c)=(x+c)[(x+w)(x+d)+1]\in\la x+c\ra  .
\end{equation}
Similarly,
\begin{equation}\label{p3}
\tut(x+d)=(x+d)[(x+w)(x+c)+1]\in\la x+d\ra  .
\end{equation}
By (\ref{p1}), (\ref{p2}) and (\ref{p3}), we obtain that $\tut(p_t(x))\in \la p_t(x)\ra  $,
and hence
$\tut$ induces an endomorphism of $k((t))[x]/\la p_t(x)\ra  $
which we continue to denote by $\tut$.
As can easily be verified, the primitive idempotents given in
(\ref{idem1}) and (\ref{idem23}) are fixed under $\tut$:
\begin{equation}\label{ei}
\tut(e_i)=e_i,\ \ i=1,2,3,
\end{equation}
whereas
\begin{equation}\label{txe}
\tut(\bx e_1)=\tut(\bx)e_1=(\bx\pi(\bx)+\bx+a)e_1=(\bx+a)e_1.
\end{equation}
Hence, $\tut$ induces an automorphism
of $K$ of order 2 while fixing the two copies of $k((t))$ pointwise.
Furthermore, one can easily verify that
$$\eta_{t=0}(\bx)=\bx^3.$$
Consequently, the automorphism
$\tut$ of $[{kC_4}]_t$ agrees with the automorphism $\eta(\btu)$ of $kC_4$ when $t=0$.
The skew polynomial ring
$$[{kC_4}]_t[y;\tut]=(k((t))[x]/\la p_t(x)\ra) [y;\tut]$$
is therefore a deformation of $kC_4[y;\eta]$.

Note that by (\ref{ei}), the
idempotents $e_i, i=1,2,3$ are central in $[{kC_4}]_t[y;\tut]$ and hence
\begin{equation}\label{sumid}
[{kC_4}]_t[y;\tut]=\bigoplus\limits_{i=1}^3[{kC_4}]_t[y;\tut]e_i.
\end{equation}

\subsection{Example}\label{example}
The following is an example for the above construction.

Put
$$a:=\frac{t+t^2+t^3}{1+t}, b:=1+t^2+t^3, c:=\frac{1}{1+t}, d:=1+t+t^2, w:=t.$$
These elements satisfy equation (\ref{assume}):
\begin{equation*}
\begin{split}
(x+w)&(x+c)(x+d)~=~
(x+t)(x+\frac{1}{1+t})(x+1+t+t^2)\\
=&~x^3+\frac{t+t^2+t^3}{1+t}x^2+(1+t^2+t^3)x+\frac{t+t^2+t^3}{1+t}
~=~x\pi(x)+a.
\end{split}
\end{equation*}
The polynomial $$\pi(x)=x^2+\frac{t+t^2+t^3}{1+t}x+{1+t^2+t^3}$$
does not admit roots in $k[[t]]/\la t^2\ra  $, thus it is irreducible over $k((t))$.

\section{A Deformation of $q(y)$}\label{hbt}
The construction of $[{kQ_8}]_t$ will be completed once the product
$\btu * \btu$ is defined.
For this purpose the polynomial $q(y)$ (\ref{q}), which determined the
ordinary multiplication $\tau^2$, will now be developed in powers of $t$.

For any non-zero element $z\in k[[t]]\setminus k[[t]]^*$, let
\begin{equation}\label{qt1}
q_t(y):=y^2+z\bx\pi(\bx)y+\bx^2+a\bx\in [kC_4]_t[y;\tut].
\end{equation}
Decomposition of (\ref{qt1}) with respect to the idempotents $e_1,e_2,e_3$ yields
\begin{equation}\label{qt}
q_t(y)=(y^2+b)e_1+[y^2+zay+c(c+a)]e_2+[y^2+zay+d(d+a)]e_3.
\end{equation}
We now show that $q_t(y)$ is in the center of $[kC_4]_t[y;\tut]:$

First, the leading term $y^2$ is central since the automorphism $\tut$ is of order 2. Next, by (\ref{ei}),
the free term $be_1+c(c+a)e_2+d(d+a)e_3$ is invariant under the action of $\tut$ and hence central.
It is left to check that the term $za(e_2+e_3)y$ is central.
Indeed, since $e_2$ and $e_3$ are $\tut$-invariant, then $za(e_2+e_3)y$ commutes both with
$[{kC_4}]_t[y;\tut]e_2$ and $[{kC_4}]_t[y;\tut]e_3$. Furthermore, by orthogonality
$$za(e_2+e_3)y\cdot[{kC_4}]_t[y;\tut]e_1=[{kC_4}]_t[y;\tut]e_1\cdot za(e_2+e_3)y=0,$$
and hence $za(e_2+e_3)y$ commutes with $[{kC_4}]_t[y;\tut]$.

Consequently, $\la q_t(y)\ra=q_t(y)[kC_4]_t[y;\tut]$ is a two-sided ideal.

Now, as can easily be deduced from (\ref{qt1}),
\begin{equation}\label{qt0}
q_{t=0}(y)=y^2+\bx^2=q(y),
\end{equation}
where the leading term $y^2$ remains unchanged.
Then $$[kQ_8]_t:=[kC_4]_t[y;\tut]/\la q_t(y)\ra$$
is a deformation of $kQ_8$,
identifying $\btu$ with $\bar{y}:=y+\la q_t(y)\ra$.

\section{Separability of $[kQ_8]_t$}\label{sep}
Finally, we need to prove that the deformed algebra $[kQ_8]_t$ is separable.
Moreover, we prove that its decomposition to simple components
over the algebraic closure of $k((t))$ resembles that of $\C Q_8$.
By (\ref{sumid}), we obtain
\begin{equation} \label{jkQ8}
[kQ_8]_t=
\bigoplus\limits_{i=1}^3[{kC_4}]_t[y;\tut]e_i/\la q_t(y)e_i \ra.
\end{equation}
We handle the three summands in (\ref{jkQ8}) separately:

By (\ref{qt}),
$$[{kC_4}]_t[y;\tut]e_1/\la q_t(y)e_1\ra\simeq K[y;\tut]/\la y^2+b \ra\simeq K^f*C_2.$$
The rightmost term is the {\it crossed product} of the group $C_2:=\la \bar{\tau}\ra$
acting faithfully on the field $K=[{kC_4}]_te_1$ via $\tut$
(\ref{txe}), with a twisting determined by the 2-cocycle $f:C_2\times C_2\to K^*$:
$$ f(1,1)=f(1,\btu)=f(\btu,1)=1,\ \, f(\btu,\btu)=b.$$
This is a central simple algebra over the subfield of invariants $k((t))$ \cite[Theorem 4.4.1]{H}.
Evidently, this simple algebra is split by $\overline{k((t))}$, i.e.
\begin{equation} \label{e1}
[{kC_4}]_t[y;\tut]e_1/\la q_t(y)e_1\ra \otimes_{k((t))} \overline{k((t))}\simeq M_2(\overline{k((t))}).
\end{equation}

Next, since $\tut$ is trivial on $[kC_4]_te_2$, the skew polynomial ring $[{kC_4}]_te_2[y;\tut]$ is actually
an ordinary polynomial ring $k((t))[y]$.
Again by (\ref{qt}),
$$[{kC_4}]_t[y;\tut]e_2/\la q_t(y)e_2\ra\simeq k((t))[y]/\la y^2+zay+c(c+a)\ra.$$
Similarly,
$$[{kC_4}]_t[y;\tut]e_3/\la q_t(y)e_3\ra\simeq k((t))[y]/\la y^2+zay+d(d+a)\ra.$$
The polynomials $y^2+zay+c(c+a)$ and $y^2+zay+d(d+a)$ are separable (since $za$ is non-zero).
Thus, both $[{kC_4}]_t[y;\tut]e_2/\la q_t(y)e_2\ra$ and $[{kC_4}]_t[y;\tut]e_3/\la q_t(y)e_3\ra$ are separable
$k((t))$-algebras, and for $i=2,3$
\begin{equation} \label{e23}
[{kC_4}]_t[y;\tut]e_i/\la q_t(y)e_i\ra \otimes_{k((t))} \overline{k((t))}
\simeq \overline{k((t))}\oplus \overline{k((t))}.
\end{equation}
Equations (\ref{jkQ8}), (\ref{e1}) and (\ref{e23}) yield
$$[kQ_8]_t\otimes_{k((t))} \overline{k((t))}\simeq \bigoplus_{i=1}^4\overline{k((t))}\oplus M_2(\overline{k((t))})$$
as required.

\section{Acknowledgement}\label{ack}
We wish to thank M. Schaps for pointing out to us that there is an error in the attempted
proof in \cite{GG} that the quaternion group is a counterexample to the DF conjecture.
Here is her explanation:
The given relations for the group algebra are incorrect. Using the notation
in pages 166-7 of \cite{GG}, if $a = 1 + i$,
$b = 1 +j$ and $z = i^2=j^2$, then $ab +ba = ij(1+z)$ while $a^2 = b^2 = 1+z$.
There is a further error later on when the matrix algebra is
deformed to four copies of the field, since a non-commutative
algebra can never have a flat deformation to a commutative algebra.

\end{document}